\documentclass[12pt,twoside]{article} 
\usepackage{amsfonts}
\usepackage{graphicx}
\usepackage[small,sc]{caption2} 

\newtheorem{theorem}{Theorem}[section]
\newtheorem{lemma}[theorem]{Lemma}
\newtheorem{proposition}[theorem]{Proposition}

\newtheorem{definition}[theorem]{Definition}

\newtheorem{rmrk}[theorem]{Remark}

\makeatletter
\@addtoreset{equation}{section}
\makeatother

\newcommand{\fig}[3] {
\medskip\smallskip
\begin{figure}[htb]
	\centering
	\includegraphics[width=#2]{#1.eps}
	\begin{minipage}[t]{0.80\linewidth} 
		\caption{#3}
		\protect\label{#1}
	\end{minipage}
\end{figure}
\medskip
}

\newenvironment{remark}
{\begin{rmrk} \em}
{\end{rmrk}}

\newcommand{\R} {{\mathbb R}}
\newcommand{\C} {{\mathbb C}}

\newcommand{\Z} {{\mathbb Z}}
\newcommand{\N} {{\mathbb N}}
\newcommand{\qed} {\hfill {\small Q.E.D.} \par\medskip}
\newcommand{\skippar} {\par\medskip}

\newcommand{\proof} {\noindent \textsc{Proof.} }
\newcommand{\proofof}[1] {\noindent \textsc{Proof of {#1}.} }
\newcommand{\article}[3] {\textsc{{#1}}, {\itshape {#2}}, {{#3}}.}
\newcommand{\book}[3] {\textsc{{#1}}, {\itshape {#2}}, {{#3}}.}
\newcommand{\vol} {\textbf}
\newcommand{\eps} {\varepsilon}
\newcommand{\rset}[2] {\left\{ #1 \: \left| \: #2 \right. \! \right\} }

\renewcommand{\iff} {if and only if\ }

\newcommand{\fn} {function}
\newcommand{\me} {measure}

\newcommand{\erg} {ergodic}
\newcommand{\sy} {system}
\renewcommand{\o} {orbit}
\newcommand{\nps} {\mathcal{N}}
\newcommand{\ph} {\varphi}
\newcommand{\pr} {probability}
\newcommand{\ra} {random}
\newcommand{\rw} {random walk}
\newcommand{\en} {environment}
\newcommand{\es} {\mathcal{E}}
\renewcommand{\P} {{\mathbb P}}
\newcommand{\bfe} {\mathbf{E}}
\newcommand{\bfn} {\mathbf{N}}
\newcommand{\bfw} {\mathbf{W}}
\newcommand{\bfs} {\mathbf{S}}
\newcommand{\bfr} {\mathbf{R}}
\newcommand{\bff} {\mathbf{F}}
\newcommand{\bfl} {\mathbf{L}}
\newcommand{\bfb} {\mathbf{B}}

\begin{document}

\title{\textbf{Recurrence for persistent random walks in two dimensions}}

\author{\textsc{Marco Lenci}
\thanks{
Department of Mathematical Sciences,
Stevens Institute of Technology, 
Hoboken, NJ 07030, U.S.A. \ 
E-mail: \texttt{mlenci@math.stevens.edu} 
} }

\date{July 2005}

\maketitle

\begin{abstract}
  We discuss the question of recurrence for persistent, or Newtonian,
  random walks in $\Z^2$, i.e., random walks whose transition
  probabilities depend both on the walker's position and incoming
  direction. We use results by T\'oth and Schmidt--Conze to prove
  recurrence for a large class of such processes, including all
  ``invertible'' walks in elliptic random environments. Furthermore,
  rewriting our Newtonian walks as ordinary random walks in a suitable
  graph, we gain a better idea of the geometric features of the
  problem, and obtain further examples of recurrence.

  \bigskip\noindent Mathematics Subject Classification: 60G50, 37B20,
  60K37, 82C41.
\end{abstract}

\section{Introduction}
\label{sec-intro}

A persistent \rw\ (PRW) in a lattice is a second order Markov chain on
that lattice. In other words, if $\{ X_n \}_{n\in\N}$ is a realization
of the chain (a \emph{walk}), the transition probabilities at time $n$
depend not only on the position $X_n$, but also on the incoming
direction $X_n - X_{n-1}$ (in a sense, the walker's ``velocity''). For
this reason, one can term this type of \rw\ `Newtonian'.

In this note we treat exclusively PRWs in $\Z^2$ with nearest-neighbor
transitions, or \emph{jumps}, but we consider both space-homogeneous
and inhomogeneous transition probababilities (in jargon, \emph{\en
s}). In the latter framework, we pay particular attention to \ra\ \en
s (REs).

Our goal is to give general ideas and specific techniques to prove
recurrence for these \rw s. In particular, using results by Schmidt
and Conze \cite{sch,co} and by T\'oth \cite{t}, we establish
recurrence for a large class or PRWRE. Also, we present an \emph{ad
hoc} construction that maps our nearest-neighbor Newtonian \rw s in
$\Z^2$ into first order \rw s in a certain graph $\Gamma$. For some
examples, this isomorphism easily entails recurrence, or lack
thereof. We hope that this construction can be fertile ground for
further study.

\skippar

We lay down the notation: A persistent \rw\ in $\Z^2$ is a Markov
chain $P^\omega$ on $\Z^2 \times \Delta$, where $\Delta := \{ \pm e_1,
\pm e_2 \}$ ($e_1, e_2$ are the generators of $\Z^2$) and $\omega$,
the \en, is an element of
\begin{equation}
  \es := M_\Delta^{\Z^2} = \rset{ \omega = \{ \omega_x \}_{x\in \Z^2}
  } { \omega_x \in M_\Delta }
  \label{def-es}
\end{equation}
Here $M_\Delta$ denotes the set of all $4 \times 4$ stochastic
matrices indexed by the elements of $\Delta$. If $\xi := (x_0,d_0) \in
\Z^2 \times \Delta$ is an initial condition, the Markov chain is
uniquely determined by the laws
\begin{eqnarray}
  && \hspace*{-30pt} P_\xi^\omega \{ (X_0, D_0) = (x_0,d_0) \} = 1;
  \label{law1} \\
  && \hspace*{-30pt} P_\xi^\omega \{ (X_{n+1}, D_{n+1})=(x',d') \,|\, 
  (X_n, D_n) = (x,d) \} = \left\{ 
  \begin{array}{rl}
    \omega_x (d,d'), & \mbox{if } x' \!=\! x \!+\! d'; \\
    0, & \mbox{otherwise.} 
  \end{array}
  \right. 
  \label{law2}
\end{eqnarray}
The \en\ and the walk are called \emph{homogeneous} if $\omega_x =$
constant.

If a probability $\Pi$ is defined on $\es$---more precisely on
$\Sigma$, the infinite-product $\sigma$-algebra naturally induced by
definition (\ref{def-es})---and $\omega$ is a \ra\ element of $(\es,
\Sigma, \Pi)$, then we have a measured family (an \emph{ensemble}) $\{
P^\omega \}_{\omega \in \es}$ of PRWs. In this case we speak of
`persistent \rw\ in a \ra\ \en'. As it is physically reasonable, we
will always assume $\Pi$ to be invariant for the natural action of
$\Z^2$ on $\es$. In fact, $\Pi$ is generally taken \erg\ for that
action; e.g., $\Pi$ can be the product of $\Z^2$ copies of a \pr\
measure on $M_\Delta$ (case of i.i.d.\ transition probabilities).

\fig{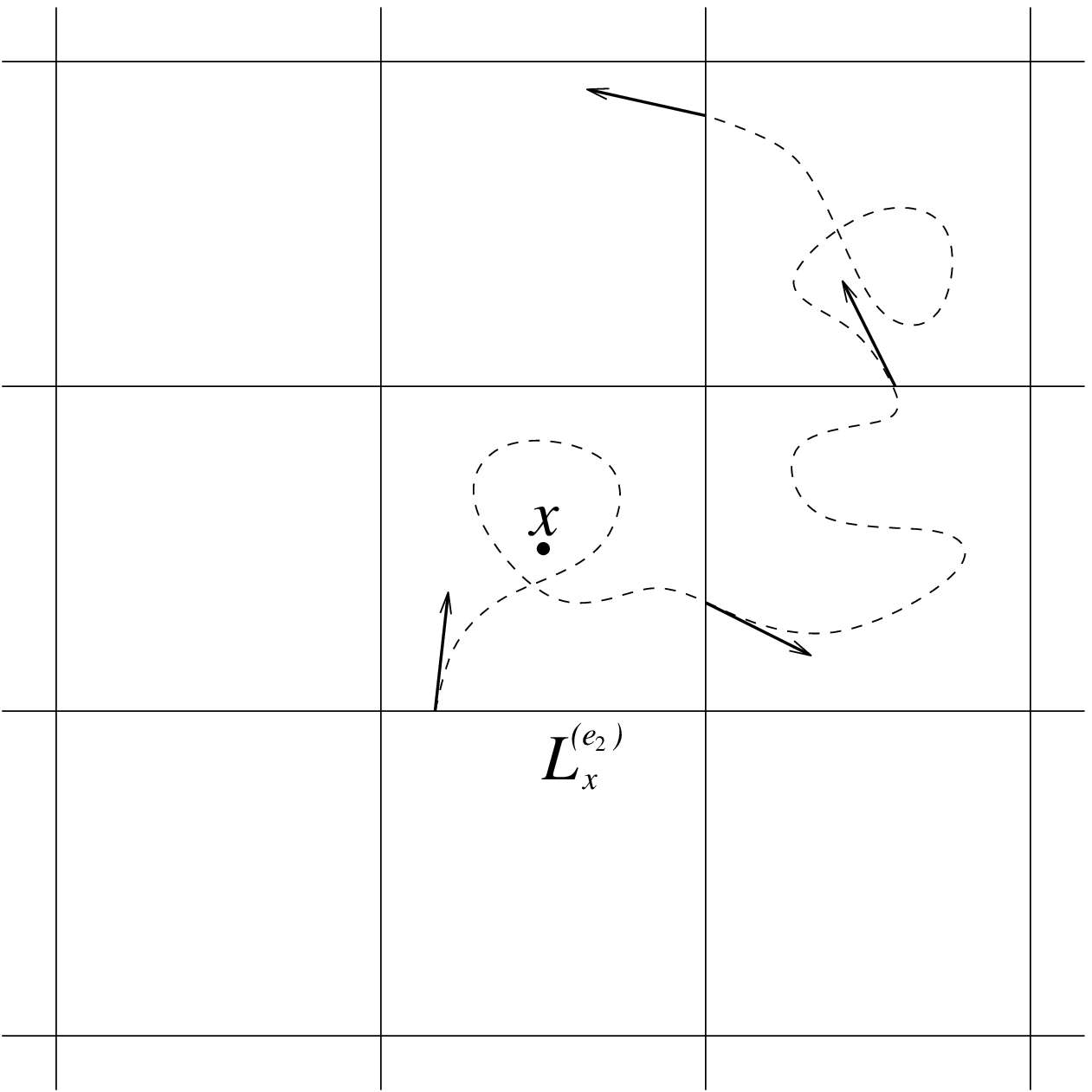}{3.2in}{Construction of the cross-section
$\nps$. $L_x^{(e_2)}$ is the bottom side of the unit square centered
in $x \in \Z^2$.}

Persistent \rw s were introduced by Sz\'asz and T\'oth in the 1980s
\cite{szt,t} as stochastic models for systems of a particle moving
chaotically in a $\nu$-dimensional space under the influence of a
certain force, for example the gradient of \ra\ potential.
Specializing to $\nu=2$ and observing Fig.~\ref{fig-rwlg-1}, call
$L_x^{(d)}$ the side of the unit square centered in $x \in \Z^2$ that
has $d$ as its inner normal. Indicating by $q,p \in \R^2$,
respectively, the position and the momentum of the particle, consider
the cross section
\begin{equation}
  \nps = \bigcup_{\scriptstyle x \in \Z^2 \atop \scriptstyle d \in
  \Delta} \nps_x^{(d)} := \bigcup_{\scriptstyle x \in \Z^2 \atop
  \scriptstyle d \in \Delta} \rset{(q,p) \in \R^2 \times \R^2} {q \in
  L_x^{(d)}, \ p \cdot d >0 }.
  \label{nps}
\end{equation}
If $T: \nps \longrightarrow \nps$ is the first-return map generated by
the dynamics of the particle, then, in general,
\begin{equation}
  (q,p) \in \nps_x^{(d)} \quad \Longrightarrow \quad T(q,p) \in
  \nps_{x+f}^{(f)}, \mbox{ for some } f \in \Delta
\end{equation}
(excluding of course those orbits that intersect $L_{x+f}^{(f)}$
tangentially or in a corner). If one assumes that, relative to some
distribution of initial conditions, the dynamics inside each unit
square is so chaotic that the particle's outgoing side appears to
depend only on the incoming side and not on the exact incoming
coordinates, then one is approximately justified in replacing a true
\o\ $\{ T^n(q,p) \}$ of the map with a realization $\{ (X_n, D_n) \}$
of the \rw.

\skippar

We call a PRW \emph{isotropic} if, $\forall x \in \Z^2$, $d' \in
\Delta$,
\begin{equation}
  \sum_{d \in \Delta} \omega_x (d,d') = 1;
  \label{isotropic}
\end{equation}
that is, the $\omega_x$ are doubly stochastic matrices. In this case,
the process can be inverted by means of another PRW. Furthermore, the
family $\{ P_\xi^\omega \}_{\xi \in \Z^2 \times \Delta}$ is
time-invariant in the following sense: For all positive integers $k$,
\begin{eqnarray}
  && \hspace*{-20pt} P_\xi^\omega \{ (X_1,D_1) = (x_1,d_1), \ldots,
  (X_n,D_n) = (x_n,d_n) \} = \\
  &&  \hspace*{-20pt} = \sum_\eta P_\eta^\omega \{ (X_k,D_k) = \xi,
  (X_{k+1},D_{k+1}) = (x_1,d_1), \ldots, (X_{k+n},D_{k+n}) =
  (x_n,d_n) \}. \nonumber
\end{eqnarray}
This type of situation simulates, for instance, a Hamiltonian \sy\
defined by a \ra\ potential $V_\omega(q)$ that is zero on every
$L_x^{(d)}$ (cf. the \ra\ Lorentz gases of \cite{l2}). In this \sy, if
one takes into account the conservation of energy (redefining
$\nps_x^{(d)}$ in (\ref{nps}) with the additional condition $|p|^2 =$
constant), one sees that the Liouville \me\ induces on $\nps$ an
invariant \me\ $\mu$ that assigns the same weight to every
$\nps_x^{(d)}$, whence the term `isotropic'. Clearly $\mu(\nps) =
\infty$.

A PRW is said to be \emph{elliptic} if there is an $\eps > 0$ such
that, $\forall x, d, d'$,
\begin{equation}
  \omega_x (d,d') \ge \eps.
  \label{elliptic}
\end{equation}
A PRWRE is elliptic if (\ref{elliptic}) holds with an $\eps$
independent of $\omega \in \es$.

\skippar

The question of recurrence for ordinary (i.e., first order) RWRE is
settled in dimension one but open and intensively pursued for $\nu \ge
2$ \cite{z}. The seminal work of Schmidt \cite{sch} and Conze
\cite{co} has shown that, in dimension two, recurrence is implied by
the Central Limit Theorem (plus \erg ity). This result has been
established only in particular cases, the most famous of which are
perhaps Lawler's `balanced' walks \cite{la} (where the transition
probabilities in the directions $v$ and $-v$ are the same). We also
mention the improvement by Komorowski and Olla \cite{ko} for certain
`two-fold stochastic' walks. It should be noted that all these results
require the ellipticity condition, the analog of (\ref{elliptic}) for
ordinary RWRE.

\skippar

This is the plan of the paper: In Section \ref{sec-prel} we lay out
some basic definitions together with a simple recurrence result for
homogeneous \en s.  In Section \ref{sec-toth} we prove that all
isotropic elliptic PRWRE (and more) are recurrent. In Section
\ref{sec-dual} we make a step further by introducing $\Gamma$, the
dual graph for nearest-neighbor PRWs in $\Z^2$. PRWs in $\Z^2$ become
first order RWs in $\Gamma$. We see that $\Gamma$ is not outrageously
complicated and, in some cases, recurrence is easy to establish.
Finally, in the Appendix, we present a recurrence result for the
Manhattan lattice which will be used (and generalized) throughout
Section \ref{sec-dual}.

\bigskip

\noindent
\textbf{Acknowledgments.} I am greatly indebted to Eric Carlen, to
whom I owe the dual graph idea (Section \ref{sec-dual}).  I also wish
to thank Massimo Campanino, Nadine Guillotin, Firas Rassoul-Agha,
Domokos Sz\'asz, and Balint T\'oth. This research was partially
supported by NSF Grant DMS-0405439.

\section{Preliminary definitions and homogeneous walks}
\label{sec-prel}

In our Newtonian RW it is necessary to distinguish between
\emph{absolute} directions in $\Z^2$ and \emph{relative} directions to
the motion of the \ra\ walker. We call the former East, North, West,
and South, and use the symbols
\begin{equation}
  \bfe := e_1, \qquad \bfn := e_2, \qquad \bfw := -e_1, \qquad \bfs :=
  -e_2.
 \label{enws}
\end{equation}
For the latter, we use the terms Right, Forward, Left, and
Backward. These are \fn s of a given direction $d \in \Delta$;
precisely, if $d^\perp \in \Delta$ is such that $(d, d^\perp)$ has the
same orientation as $(e_1, e_2)$, they are defined as
\begin{equation}
  \bfr(d) := -d^\perp, \qquad \bff(d) := d, \qquad \bfl(d) := d^\perp,
  \qquad \bfb(d) := -d.
 \label{rflb}
\end{equation}
In ascending order of strength, the following are the properties that
we are interested in, for our \rw s:

\begin{definition}
  The PRW in the \en\ $\omega \in \es$ is said to have
  \textbf{asymptotic velocity} $v \in \R^2$ if, $\forall \xi \in \Z^2
  \times \Delta$,
  \begin{displaymath}
    \lim_{n \to +\infty} \frac{X_n}n = v \qquad 
    P_\xi^\omega\mbox{\rm -almost surely}.
  \end{displaymath}
  In case $v \ne 0$, one says that the walk is \textbf{ballistic}.
  \label{def-as-vel}
\end{definition}

\begin{definition}
  The PRW in $\omega$ is \textbf{diffusive} if there exist two
  constants $C_2 > C_1 > 0$ such that, $\forall \xi \in \Z^2 \times
  \Delta$,
  \begin{displaymath}
    C_1 \, n \le E_\xi^\omega \left( |X_n|^2 \right) \le C_2 \, n,
  \end{displaymath}
  where $E_\xi^\omega$ is the expectation w.r.t.\ $P_\xi^\omega$.
  If, in addition, $X_n / \sqrt{n}$ tends in distribution to a
  non-degenerate centered Gaussian, then we say that the walk
  satisfies the \textbf{Central Limit Theorem (CLT)} with zero 
  asymptotic velocity.
  \label{def-diff}
\end{definition}

\begin{definition}
  The PRW in $\omega$ is \textbf{recurrent} if, $\forall \xi \in \Z^2
  \times \Delta$,
  \begin{displaymath}
    P_\xi^\omega \{ (X_n, D_n) = \xi = (X_0, D_0), \mbox{ \rm for 
    some } n>0 \} = 1.
  \end{displaymath}
  By the Markov property, this clearly implies that, almost surely,
  $(X_n, D_n) = (X_0, D_0)$ for infinitely many $n$. On the other
  hand, the walk is \textbf{transient} if, for every finite $B \subset
  \Z^2 \times \Delta$, there exists almost surely an $m>0$ such that
  $(X_n, D_n) \not\in B$ $\forall n\ge m$.
  \label{def-rec}
\end{definition}

\begin{remark}
  In the case of \ra\ \en\ $(\es, \Sigma, \Pi)$, all these definitions
  are meant for $\Pi$-a.e.\ $\omega \in \es$.
  \label{rk-ae}
\end{remark}

In the homogeneous case, that is, when $\omega_x =$ constant, these
properties are rather easily investigated with the methods of finite
Markov chains. Nevertheless, for the sake of completeness, we present
here a rather general result for that case.

\begin{theorem}
  Assume that a homogeneous PRW, (i.e., $\omega_x = \omega_0$ $\forall
  x \in \Z^2$) is defined by an irreducible aperiodic stochastic
  matrix $\omega_0 \in M_\Delta$. A necessary and sufficient condition
  for the walk to satisfy the CLT with zero asymptotic velocity, and
  be recurrent, is that
  \begin{eqnarray}
    &\phantom{=}& \frac{1 - \omega_0(\bfe,\bfe) - \omega_0(\bfw,\bfe)}
    {\omega_0(\bfn,\bfe) + \omega_0(\bfs,\bfe)} = \frac{1 -
    \omega_0(\bfe,\bfw) - \omega_0(\bfw,\bfw)} {\omega_0(\bfn,\bfw) +
    \omega_0(\bfs,\bfw)} = \nonumber \\
    &=& \frac{\omega_0(\bfe,\bfn) + \omega_0(\bfw,\bfn)} {1 -
    \omega_0(\bfn,\bfn) - \omega_0(\bfs,\bfn)} =
    \frac{\omega_0(\bfe,\bfs) + \omega_0(\bfw,\bfs)} {1 -
    \omega_0(\bfn,\bfs) - \omega_0(\bfs,\bfs)}.  \label{hom-conds}
  \end{eqnarray}
  If the walk is isotropic, as per definition (\ref{isotropic}), the
  above condition is clearly verified.
  \label{thm-hom}
\end{theorem}

\proof As the \en\ is homogeneous, fix without loss of generality $\xi
= (x_0, d_0) = (0, e_1)$. The statistical properties of
\begin{equation}
  X_n = \sum_{k=0}^{n-1} D_k
  \label{xn}
\end{equation}
are simply inferred from the statistical properties of the finite
Markov chain $\{ D_n \}$, governed by the stochastic matrix
$\omega_0$. It is well known that a finite, irreducible, aperiodic
Markov chain satisfies the CLT for the number of visits in each
state. More precisely, if $\pi = (\pi^\bfe, \pi^\bfn, \pi^\bfw,
\pi^\bfs)$ is the (unique) stationary vector of $\omega_0$, determined
by
\begin{equation}
  \sum_{d \in \Delta} \pi^d \, \omega_0(d,d') = \pi^{d'} \qquad
  \forall d' \in \Delta,
  \label{thm-hom-1}
\end{equation}
and $C_n = (C_n^\bfe, C_n^\bfn, C_n^\bfw, C_n^\bfs)$ is the so-called
\emph{counting vector}, given by
\begin{equation}
  C_n^d := \# \rset{0 \le k \le n-1} {D_k = d};
\end{equation}
then
\begin{equation}
  \frac{ C_n - n \pi } {\sqrt{n}} 
\end{equation}
converges in distribution, as $n \to +\infty$, to a multivariate
centered Gaussian \cite{rs}. By linearity, then, $(X_n - n
V)/\sqrt{n}$, with
\begin{equation}
  V := \pi^\bfe \, \bfe + \pi^\bfn \, \bfn + \pi^\bfw \, \bfw +
  \pi^\bfs \, \bfs,
\end{equation}
also converges to a centered Gaussian. Therefore, the RW will satisfy
the CLT as per Definition \ref{def-diff} \iff $V=0$, that is, \iff
$\pi^\bfe = \pi^\bfw$ and $\pi^\bfn = \pi^\bfs$.

Let us prove the sufficiency of condition (\ref{hom-conds}).  Denote
by $\lambda \in [0,+\infty]$ the value of the four equal expressions
in (\ref{hom-conds}). For each such $\lambda$, there is a unique pair
$p_1, p_2 \ge 0$ such that
\begin{equation}
  p_1 + p_2 = \frac12, \qquad \frac{p_2}{p_1} = \lambda.
  \label{thm-hom-2}
\end{equation}
Is it easy to check that (\ref{hom-conds}) is equivalent to $p :=
(p_1, p_2, p_1, p_2)$ verifying (\ref{thm-hom-1}), and thus being
equal to $\pi$.  This last sentence also shows the necessity of
(\ref{hom-conds}) if $V=0$, that is, if $\pi$ is of the form $(p_1,
p_2, p_1, p_2)$.  

As for the recurrence, let us reinterpret our Markov chain as the
dynamical \sy\ $(\Delta^\N, \sigma, P)$, where $\sigma$ is the left
shift on $\Delta^\N$ and $P$ is the Markov \pr\ induced on $\Delta^\N$
by the stochastic matrix $\omega_0$ and the stochastic vector
$\pi$. It is well known \cite{rs} that if $\omega_0$ is irreducible
and aperiodic (a.k.a.\ \erg) then this \sy\ is \erg. Furthermore, by
the above, the \emph{cocycle} (\ref{xn}) satisfies the CLT with zero
mean. Hence, by Conze's theorem on the recurrence of dynamical
cocycles \cite{co}, $X_n = 0$, $P$-a.s., for infinitely many $n$. For
at least one (and thus infinitely many) such $n$, it must also be that
$D_n = e_1$ (by the Markov property and the fact that, if $(X_n, D_n)
= (0,d)$, there is a positive \pr\ that $(X_k, D_k) = (0,e_1)$ for
some $k>n$).
\qed

\section{T\'oth environments}
\label{sec-toth}

A large class of PRWRE can be proven recurrent by combining the
Schmidt--Conze theorem with a previous result by T\'oth \cite{t}. The
latter extends Kipnis and Varadhan's proof of the invariant principle
\cite{kv} (namely, convergence to Brownian motion in the proper
scaling) to isotropic elliptic PRWRE, and others.

To properly state the theorem, we need a couple of definitions: For $y
\in \Z^2$ and $\omega = \{ \omega_x \}_{x \in \Z^2}$,
\begin{equation}
  \left( \tau_y \, \omega \right)_x := \omega_{x+y}
\end{equation}
represents the natural action of $\Z^2$ on $\es$. Also, let $W$ denote
the transition matrix of the standard \rw, namely, 
\begin{equation}
  W(d,d') = \frac14 \qquad \forall d,d' \in \Delta.
\end{equation}

\begin{theorem}
  Assume that $\{ \tau_y \}$ is an \erg\ group of automorphisms on
  $(\es, \Sigma, \Pi)$. Suppose also that, $\Pi$-a.s., $\omega$ is
  isotropic and,
  \begin{equation}
    \forall x \in \Z^2, \qquad \| \omega_x - W \| \le 
    1-\eps,
    \label{toth-cond-b}
  \end{equation}
  for some $\eps>0$ ($\| \,\cdot\, \|$ is the standard operator norm
  for $4 \times 4$ matrices). Then the associated PRWRE verifies the
  CLT and is recurrent.
  \label{thm-toth}
\end{theorem}

\proof First of all, due to the translation-invariance of $\Pi$, we
can assume that $X_0 = 0$. Secondly, we incorporate the randomness of
the \en\ into a unique Markov chain that describes all walks in all
\en s. This Markov chain is defined on the state space $\Delta \times
\es$ and its realizations are denoted by $\{ (D_n, \Omega_n) \}$. It
is uniquely determined by the laws
\begin{eqnarray}
  && \hspace*{-24pt} \P \{ (D_0, \Omega_0) \in A \times B \} =
  \frac{|A|}4 \: \Pi(B);
  \label{llaw1} \\
  && \hspace*{-24pt} \P \{ (D_{n+1}, \Omega_{n+1})=(d',\omega') \,|\, 
  (D_n, \Omega_n) = (d,\omega) \} = \left\{
  \begin{array}{rl}
    \omega_0 (d,d'), & \mbox{if }\ \omega' = \tau_{d'} \, \omega; \\
    0, & \mbox{otherwise.} 
  \end{array}
  \right. 
  \label{llaw2}
\end{eqnarray}
In (\ref{llaw1}), $|A|$ is the cardinality of $A \subseteq \Delta$, and
$B \in \Sigma$. In the theory of RWRE, this process is sometimes
called `the point of view of the particle'.

As in the proof of Theorem \ref{thm-hom}, we exploit the associated
dynamical \sy
\begin{equation}
  \left( [\Delta \times \es]^\N, \bar{\sigma}, \P \right), 
  \label{ds-pvp}
\end{equation}
where $\bar{\sigma}$ is the left shift on $[\Delta \times \es]^\N$.
Using (\ref{isotropic}) and the translation-invariance of $\Pi$ in
(\ref{llaw2}), we see that $(| \cdot | \otimes \Pi)/4$ (where $| \cdot
|$ is the counting \me\ on $\Delta$) is a stationary \me\ on $\Delta
\times \es$. Since this is the same as the initial distribution chosen
in (\ref{llaw1}), we conclude that $\P$ is a stationary \me\ for the
entire process, which is equivalent to saying that (\ref{ds-pvp}) is a
\me-preserving (non-invertible) dynamical \sy.

\begin{remark}
  The non-invertibility is merely our choice. Since (\ref{isotropic})
  implies that the process can be defined for negative times too,
  using the transposed environment, one could have introduced an
  invertible dynamics on $[\Delta \times \es]^\Z$, also preserving
  $\P$. In this case, Schmidt's result on recurrent cocycles
  \cite{sch} would apply as well as Conze's \cite{co}.
\end{remark}

In \cite{t} it is proved that this dynamical \sy\ is \erg\ and, more
importantly, that the associated RW satisfies the invariance
principle. This means that, defining $R_n(0) := 0$ and, for $j = 1, 2,
\ldots, n$,
\begin{equation}
  R_n \! \left( \frac{j}{n} \right) := \frac1 {\sqrt{n}}
  \sum_{k=0}^{j-1} D_k,
\end{equation}
and extending the definition to all $t \in [0,1]$ by means of linear
interpolation, the distribution of the \ra\ path $R_n: [0,1]
\longrightarrow \R^2$ converges weakly to that of a Brownian
motion. That is, if $f$ is a continuous \fn\ of the space
$\mathcal{C}([0,1]; \R^2)$ of all such paths (endowed with the usual
weak* topology), then
\begin{equation}
  \int_{[\Delta \times \es]^\N} f(R_n) \, d\P -
  \int_{\mathcal{C}([0,1]; \R^2)} f \, d\mathcal{W}_S \:
  \longrightarrow \: 0,
  \label{inv-pr}
\end{equation}
as $n \to +\infty$. Here $\mathcal{W}_S$ is the 2D Wiener \me\ with
zero drift and variance $S$ (by isotropy $S$ must be a positive
multiple of the identity). 

It is then a trivial corollary that $X_n$, which is a cocycle of
(\ref{ds-pvp}) by means of (\ref{xn}), satisfies the CLT. This and the
\erg ity of (\ref{ds-pvp}) imply, via \cite{co}, that $X_n = 0$ for
infinitely many $n$, $\P$-a.s. 

We need to prove that also $(X_n, D_n) = (0, D_0)$, for infinitely
many $n$, $\P$-a.s.  For a.e.\ $\{ (D_n, \Omega_n) \} \in [\Delta
\times \es]^\N$ it makes sense to denote by $k$ the smallest $n>0$
such that $X_n = 0$, and set
\begin{equation}
  \hat{\sigma} \{ (D_n, \Omega_n) \} := \bar{\sigma}^k \{ (D_n,
  \Omega_n) \}.
\end{equation}
The map $\hat{\sigma}$ preserves $\P$ because $\bar{\sigma}$
does. Therefore, by the Poincar\'e Theorem applied to $( [\Delta
\times \es]^\N, \hat{\sigma}, \P )$, a point that starts in $\{ D_0 =
d_0 \}$ will return there at infinitely many times $k$ with $X_k = 0$,
yielding the sought assertion. 

Which, as it is evident, is equivalent to the saying that, for every
$\xi = (0,d_0)$ and $\Pi$-a.e.~$\omega \in \es$, $\{ (X_n, D_n) \}$ is
recurrent w.r.t.~$P_\xi^\omega$.  
\qed

The next lemma shows that the strongest requirement in condition
(\ref{toth-cond-b}) is not the strictness of the bound, but its
uniformity.

\begin{lemma}
  If $Q, Q^T \in M_\Delta$ then $\| Q - W \| \le 1$. The inequality is
  strict \iff $Q^T Q$ is irreducible and aperiodic.
  \label{lem-cond-b}
\end{lemma}

\begin{remark}
  The irreducibility and aperiodicity of $Q^T Q$ are implied by a
  number of conditions given directly on $Q$:
  \begin{enumerate}
    \item {\itshape All the entries of $Q$ are positive} (which, for
    want of a better term, we may call `local ellipticity'). As a
    matter of fact, it is enough that {\itshape no column of $Q$
    contains more than 1 zero}. This implies that the entries of
    $Q^T Q$ are positive.

    \item {\itshape $Q$ is irreducible, aperiodic and normal}. In this
    case, the entries of $(Q^T Q)^k = (Q^T)^k Q^k$ are positive,
    for $k$ large enough.

    \item {\itshape $Q$ is irreducible, aperiodic, and its diagonal
    entries are positive}. In fact, in this case, $Q$ can be rewritten
    as $Q = \rho \mathbf{1} + \bar{Q}$, where $\rho>0$, $\mathbf{1}$
    is the identity, and all the entries of $\bar{Q}$ are
    non-negative.  Therefore, the entries of
    \begin{equation}
      (Q^T Q)^k = (\rho \mathbf{1} + \bar{Q}^T) Q \,\cdots\, (\rho
      \mathbf{1} + \bar{Q}^T) Q = \rho^k Q^k + \,\cdots 
    \end{equation}
    are positive, for large $k$ (the remainder terms all have
    non-negative entries).
  \end{enumerate}
  It would be a mistake, however, to conjecture that for all doubly
  stochastic, irreducible, aperiodic $Q$, $\| Q - W \| < 1$. An
  exception, for example, is
  \begin{equation}
    \left[
    \begin{array}{cccc}
      0 & 0 & 0 & 1 \\
      0 & \frac12 & \frac12 & 0 \\
      1 & 0 & 0 & 0 \\
      0 & \frac12 & \frac12 & 0
    \end{array}
    \right].
  \end{equation}
  \label{rk-cond-b}
\end{remark}

\proofof{Lemma \ref{lem-cond-b}} We start by noting that $W$ is an
orthogonal projection in $\C^4$.  The hypothesis that $Q$ is doubly
stochastic is equivalent to
\begin{equation}
  Q W = W Q = W.
\end{equation}
This implies that the eigenspaces of $W$, namely, $E := \mathrm{span}
\{ (1,1,1,1) \}$ and $E^\perp$, are invariant w.r.t.~$Q$ and
$Q^T$. The action of $Q - W$, relative to the decomposition $\C^4 = E
\oplus E^\perp$, is easy to calculate, using that $Q$ is stochastic:
\begin{equation}
  Q - W = (Q - 1)|_E \oplus Q|_{E^\perp} =  Q|_{E^\perp}.
\end{equation}
Now
\begin{equation}
  \| Q|_{E^\perp} \|^2 = \| (Q^T Q)|_{E^\perp} \| \le 1
  \label{l-cond-b-2}
\end{equation}
by the Perron-Frobenius Theorem applied to the self-adjoint stochastic
matrix $Q^T Q$. By the same theorem, $Q^T Q$ is irreducible and
aperiodic \iff its eigenvalues for the subspace $E^\perp$ have all
modulus $<1$, i.e., \iff $ \| (Q^T Q)|_{E^\perp} \| < 1$.
\qed

As an illustration of the strength and limitations of Theorem
\ref{thm-toth}, let us consider a few examples. Clearly, the theorem
applies to all isotropic and elliptic \ra\ \en s.

Moving on, for all $x \in \Z^2$, produce i.i.d.\ \ra\ vectors
$(\zeta_x, \zeta'_x, \zeta''_x)$ with $\zeta_x + \zeta'_x + \zeta''_x
= 1$. Assume nothing about the \pr\ law for the stochastic vector,
except that its components must be larger than a certain
$\eps>0$. Then set
\begin{equation}
  \begin{array}{rlllll}
    (\forall d \in \Delta) && \omega_x(d, \bff) = \zeta_x, &
    \omega_x(d, \bfl) = \zeta'_x, & \omega_x(d, \bfr) = \zeta''_x, &
    \omega_x(d, \bfb) = 0.
  \end{array}
  \label{flr-env}
\end{equation}
where we have abbreviated notation (\ref{rflb}) as $\bff = \bff(d)$,
$\bfl = \bfl(d)$, etc. This construction defines a RE $(\es, \Sigma,
\Pi)$ that is evidently isotropic and verifies the hypotheses of
Theorem \ref{thm-toth}---condition (\ref{toth-cond-b}) holds via
Remark \ref{rk-cond-b}, pt.~1 (or 3).

The same holds if the roles of $\bff$ and $\bfb$ are swapped in
(\ref{flr-env}).

On the other hand, if, $\Pi$-a.s., $\omega_x (d,\pm d)=0$ $\forall
x,d$, then Theorem \ref{thm-toth} does not apply because $\omega_x^T
\omega_x$ is reducible and $\| \omega_s - W \| = 1$. We may call any
such \en\ `Left-Right \en' (more about this in Section
\ref{sec-dual}).

\skippar

Even conceding on this last example and similar ones, the real
shortcomings of Theorem \ref{thm-toth} are two---or maybe three.
First, the uniformity of assumption (\ref{toth-cond-b}), which seems
too much of an absolute requirement for a probabilistic problem (on
the other hand, in the realm of ordinary RWs in 2D, this is mirrored
by the ellipticity condition, which is an hypothesis of all general
recurrence results known to this day \cite{z}). Second, and more
important, Theorem \ref{thm-toth} does not cover non-isotropic
REs. Third, and least relevant, it implies nothing about a given,
specific, inhomogeneous \en.

In the next section we improve the situation somewhat.

\section{The dual graph}
\label{sec-dual}

Suppose that we want to visualize a persistent RW in $\Z^2$ as a first
order RW in $\Z^2 \times \Delta$.  We quickly realize that the most
immediate mental image of $\Z^2 \times \Delta$, that is, four stacked
copies of $\Z^2$, is not really the best to work with because, due to
the costraints (\ref{law2}), from a given $(x,d) \in \Z^2 \times
\Delta$ the walker can only transit to the four other sites
$(x+d',d')$. Once we draw, for every $x, d, d'$, an oriented edge from
$(x,d)$ to $(x+d',d')$, we are left with a rather labyrinthic picture.

A better choice, it seems, is to consider $\Gamma$, the dual of the
graph $G$ whose vertices are the points of $\Z^2$ and whose oriented
edges are the pairs of nearest-neighbor relations. See
Fig.~\ref{fig-rwlg-2} for the construction of $\Gamma$.

\fig{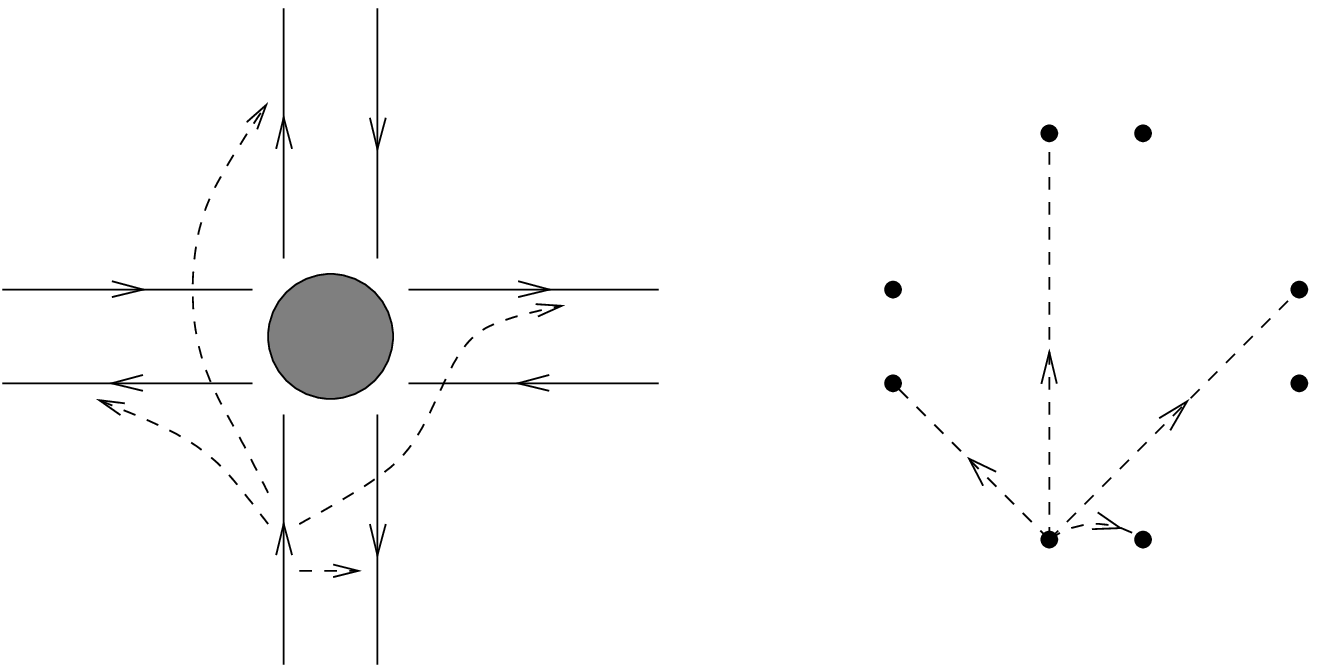}{3.7in}{Construction of $\Gamma$: The nearest-neighbor
edges of $\Z^2$ are turned into the vertices of $\Gamma$, and the
relations of possible consecutiveness between those edges are turned
into the (oriented) edges of $\Gamma$. Only one elementary building
block of $\Gamma$ is shown here; see Fig.~\ref{fig-rwlg-3} for a wider
view.}

\fig{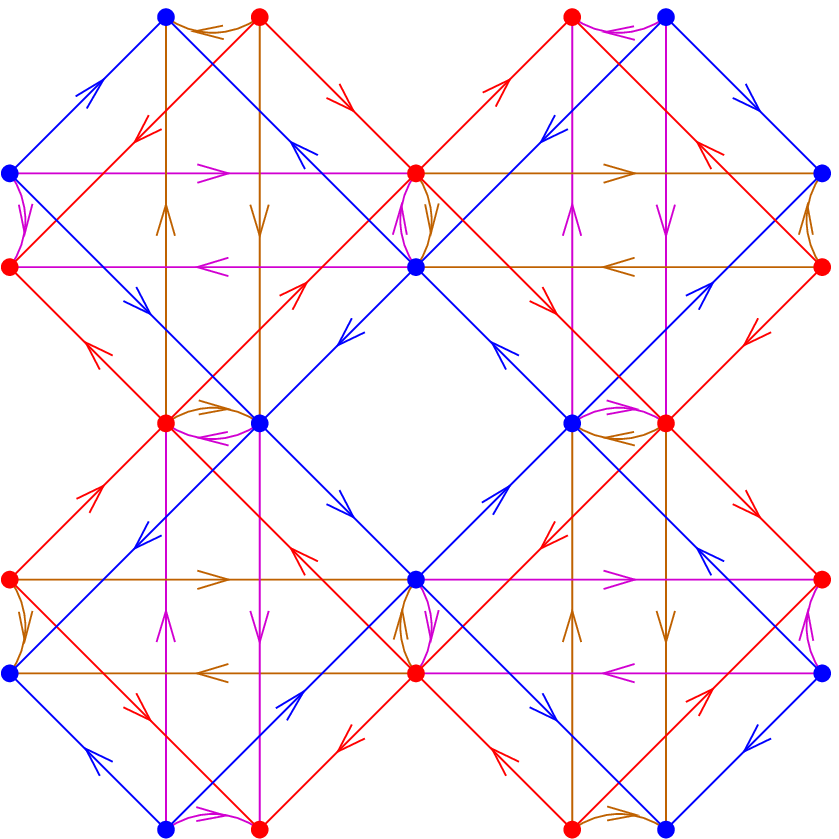}{2.8in}{A larger portion of $\Gamma$. The two embedded
Manhattan lattices $M_1, M_2$ are colored blue and red,
respectively. Edges that lead from $M_1$ to $M_2$ (that is, from blue
to red), denoted $L_+$ in the text, are colored in magenta; edges that
lead from $M_2$ to $M_1$ (red to blue), denoted $L_-$, are colored in
brown.}

If we consider the picture of $\Gamma$ given in Fig.~\ref{fig-rwlg-3},
we can make the following observations:

\begin{enumerate}

\item The diagonal edges of $\Gamma$ (those that are followed when the
particle goes Left or Right) form two connected subgraphs $M_1$ and
$M_2$, each isomorphic to the Manhattan lattice ${\mathbb M}$ (see
Appendix).

\item Calling $B$ the set of all short horizontal or vertical edges
(corresponding to the particle going Backward) and $F$ the set of all
long horizontal or vertical edges (particle going Forward), we see
that each edge in $B \cup F$ connects $M_1$ to $M_2$ or
viceversa. That is, $B \cup F = L_+ \cup L_-$, where $L_+$ is the set
of all edges leading from $M_1$ to $M_2$ and $L_-$ is the set of all
edges leading from $M_2$ to $M_1$.

\item Although $\Gamma$ is invariant for the action of $\Z^2$, the
partition $M_1 \cup M_2 \cup L_+ \cup L_-$ of its edges is only
invariant for the action of
\begin{equation}
  \Z^2_\mathrm{even} := \rset{(x^1,x^2) \in \Z^2} {x^1 + x^2 \in 2\Z}.
  \label{z2even}
\end{equation}
Actually, for $x \in \Z^2_\mathrm{odd} := \Z^2 \setminus
\Z^2_\mathrm{even}$, 
\begin{equation}
  \upsilon_x (M_1) = M_2, \quad \upsilon_x (M_2) = M_1, \quad
  \upsilon_x (L_+) = L_-, \quad \upsilon_x (L_-) = L_+, 
  \label{upsilon}
\end{equation}
where $\{ \upsilon_x \}$ is the natural action of $\Z^2$ on $\Gamma$.
This shows that $\Z^2_\mathrm{even}$ is the most complete symmetry for
the problem.

\end{enumerate}

The above considerations suggest that $\Gamma$ is better viewed as a
three-dimensional object in which the two Manhattan lattices lie one
on top of the other, with opposite orientations, as in
Fig.~\ref{fig-rwlg-4}.

\fig{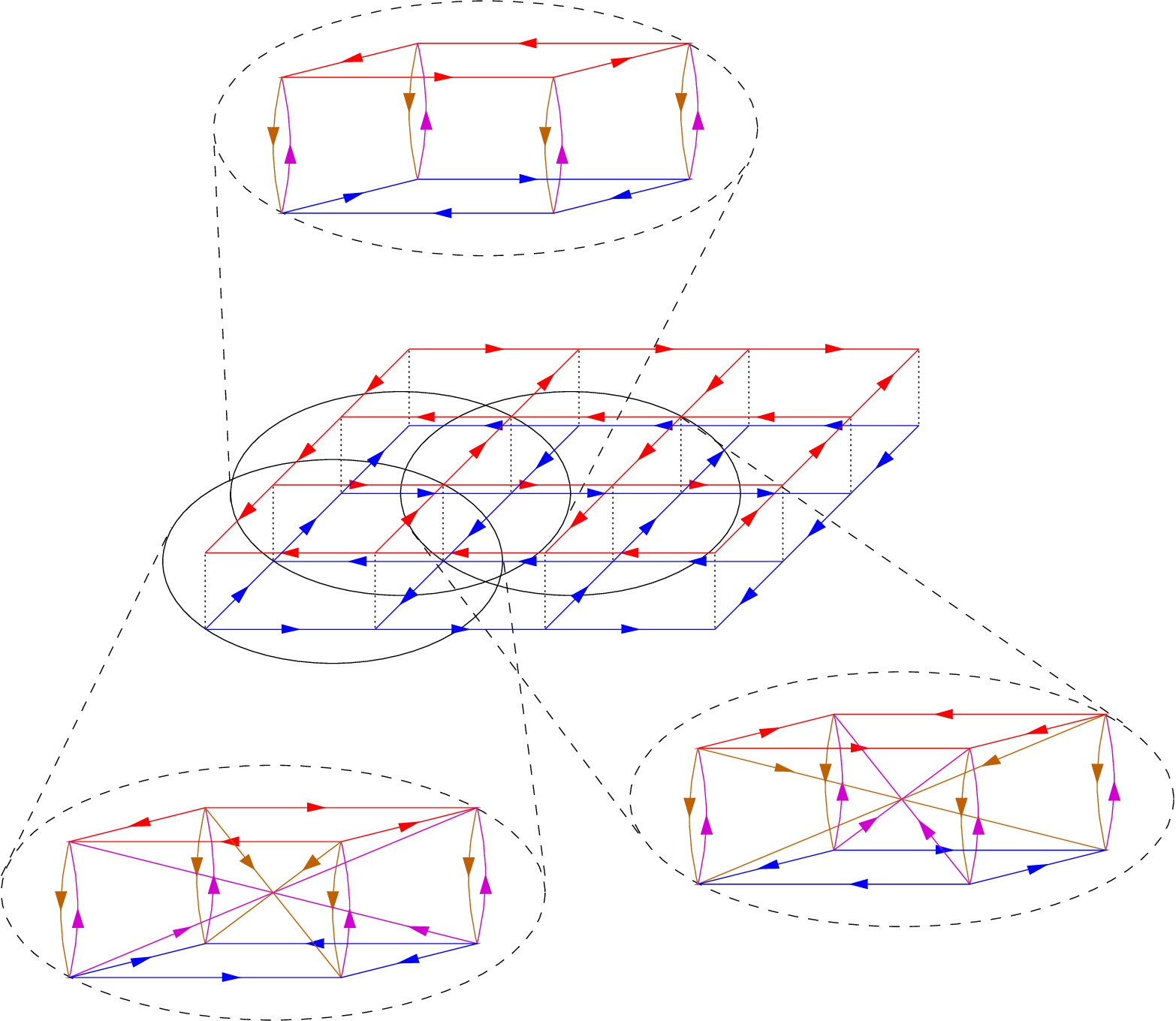}{4.6in}{Three-dimensional rendering of $\Gamma$. Here
the graph is embedded in the $(z^1, z^2, z^3)$-space and its vertices
are the points of $\Z^2 \times \{ 0,1 \}$. $M_1$ and $M_2$ lie in the
planes $\{ z^3 = 0 \}$ and $\{ z^3 = 1 \}$, respectively. Note that
all ``cubes'' do not have the same structure in terms of $L_+$ and
$L_-$ (see blowups). Indeed this structure is $(2\Z)^2$-periodic in
the displayed rendering (corresponding to
$\Z^2_\mathrm{even}$-periodicity for the problem).}

One notices that, removing $B$ and $F$ from $\Gamma$ or, equivalently,
setting $\omega_x (d, \pm d) = 0$ $\forall x,d$, splits $\Z^2 \times
\Delta$ into two parts that are never connected by an \o\ of the
RW. These are the Left-Right \en s mentioned at the end of Section
\ref{sec-toth}.

If furthermore $\omega_x (d, \pm d^\perp) = 1/2$ $\forall x,d$ (making
the \en\ homogeneous), the walk on either component is isomorphic to a
standard walk in the Manhattan lattice (see Appendix). Thus the whole
system, which we call `symmetric Left-Right RW', is recurrent. This is
a well-known fact and can be proved in several different ways (cf.\
also \cite{l2}).

One can construct \emph{ad hoc} inhomogeneous PRWs whose recurrence
properties are quickly obtained from these elementary results.

\subsection{Inhomogeneous Forward probability}

Let $\{ \zeta_x \}_{x \in \Z^2}$ be a collection of numbers from
$[0,1]$. For $x \in \Z^2_\mathrm{even}$ define,
\begin{equation}
  \begin{array}{rllll}
    (\mbox{if } d = \bfn) && \omega_x(d, \bff) = \zeta_x, & 
    \omega_x(d, \bfl) = \omega_x(d, \bfr) = (1-\zeta_x)/2, 
    & \omega_x(d, \bfb) = 0; \\
    (\mbox{if } d \ne \bfn) && \omega_x(d, \bff) = 0, & 
    \omega_x(d, \bfl) = \omega_x(d, \bfr) = 1/2, 
    & \omega_x(d, \bfb) = 0.
  \end{array}
  \label{hyp-4-1-a}
\end{equation}
(Once again, the dependence of $\bfr, \bff, \bfl, \bfb$ on $d$, as per
(\ref{rflb}), has been omitted.)  Similarly, for $x \in
\Z^2_\mathrm{odd}$,
\begin{equation}
  \begin{array}{rllll}
    (\mbox{if } d = \bfs) && \omega_x(d, \bff) = \zeta_x, & 
    \omega_x(d, \bfl) = \omega_x(d, \bfr) = (1-\zeta_x)/2, 
    & \omega_x(d, \bfb) = 0; \\
    (\mbox{if } d \ne \bfs) && \omega_x(d, \bff) = 0, & 
    \omega_x(d, \bfl) = \omega_x(d, \bfr) = 1/2, 
    & \omega_x(d, \bfb) = 0.
  \end{array}
  \label{hyp-4-1-b}
\end{equation}
In other words, this is a modification of the symmetric Left-Right
walk in the following sense: When the particle reaches a point of
$\Z^2_\mathrm{even}$ from the South, it has some \pr\ to proceed
Forward, that is, go North; when it reaches a point of
$\Z^2_\mathrm{odd}$ from the North, it has some \pr\ to go South.

Notice that these extra Forward displacements need not be
\emph{statistically balanced}. This means that, defining the
\emph{local drift} as
\begin{equation}
  \delta_\omega (x,d) := \sum_{d' \in \Delta} \omega_x(d,d') \, d',
  \label{l-drift}
\end{equation}
its $\Z^2$-average
\begin{equation}
  \lim_{\Lambda \nearrow Z^2} \, \frac1{4|\Lambda|} \sum_{\scriptstyle
  x \in \Lambda \atop \scriptstyle d \in \Delta} \delta_\omega (x,d)
  \label{avg-l-drift}
\end{equation}
may be different from zero. (Here $\Lambda$ is, e.g., a square whose
side length is going to $+\infty$, and $|\Lambda|$ is the cardinality
of $\Lambda$). Nonetheless, one can prove the following

\begin{proposition}
  If, for at least one $x \in \Z^2_\mathrm{even}$, $\zeta_x > 0$ and,
  for at least one $y \in \Z^2_\mathrm{odd}$, $\zeta_y > 0$, the PRW
  in the \en\ defined above is recurrent for all initial conditions.
  \label{prop-4-1}
\end{proposition}

\proof In this proof and for the remainder of the section we use the
3D rendering of $\Gamma$ shown in Fig.~\ref{fig-rwlg-4}. Namely, we
think of $\Gamma$ as embedded in $\R^3$ and having its vertices in the
points of $\Z^2 \times \{ 0,1 \}$. Thus, for example, $M_i$ ($i=1,2$)
is embedded in $\R^2 \times \{ i-1 \}$, $L_+$ is the set of all edges
that ``point upwards'', $B$ is the set of all vertical edges, and so
on. Notice also that the action of $\upsilon_x$, as introduced in
(\ref{upsilon}), corresponds to a translation by the vector $(x^1+x^2,
-x^1+x^2, 0)$, where $(x^1,x^2) = x$.

Studying assumptions (\ref{hyp-4-1-a})-(\ref{hyp-4-1-b}) with an eye
on Fig.~\ref{fig-rwlg-4}, we see that the only non-horizontal edges
that count (i.e., that can be traveled on with positive \pr) are those
parallel to the vector $(1,1,1)$. We call them \emph{effective} edges
of $L_+ \cup L_-$. All others are as good as removed.

Denoting by $(z^1, z^2, z^3)$ the coordinate \sy\ in $\R^3$, define
the \fn\ $\psi : \Z^2 \times \{ 0,1 \} \longrightarrow \Z^2$ as
\begin{equation}
  \psi(z^1, z^2, z^3) := (z^1 - z^3, z^2 - z^3).
  \label{proj-psi}
\end{equation}
If $\{ Z_n \}_{n \in \N}$ is a realization of our RW in $\Gamma$, with
$Z_0 = (0,0,0)$, set $Y_n := \psi(Z_n)$. Geometrically this operation
amounts to shifting $M_2$ by the vector $(-1,-1,0)$ (thus aligning the
orientations of $M_2$ and $M_1$) and projecting the walk vertically
onto $M_1$. See Fig.~\ref{fig-rwlg-5}.

\fig{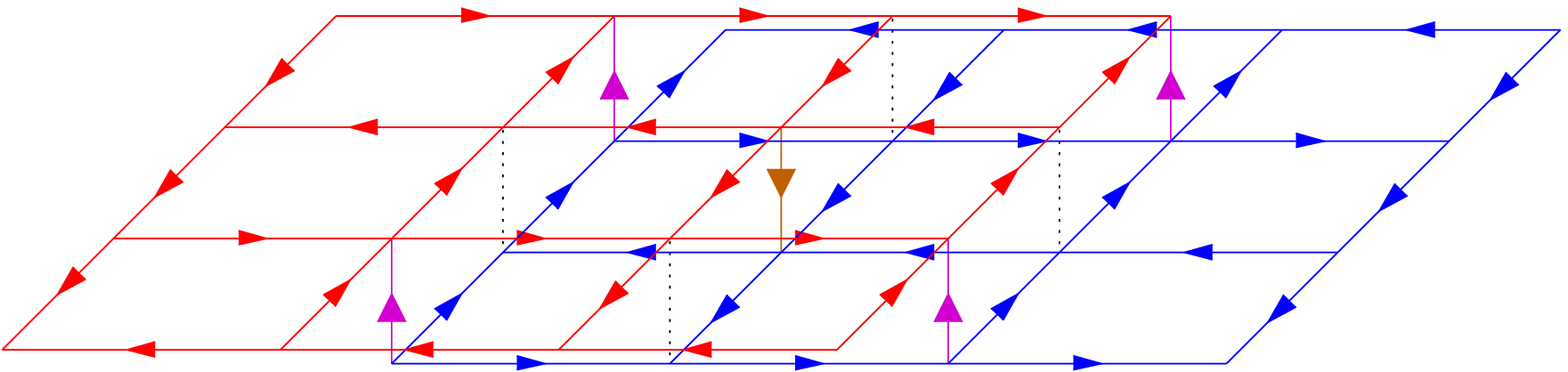}{4in}{Illustration of the proof of Proposition
\ref{prop-4-1}.}

It is clear that $\{ Y_n \}$ can only have same-site or
nearest-neighbor jumps, the former case happening when, and only when,
$Z_n$ shifts level. Define $n_0 := 0$ and, recursively,
\begin{equation}
  n_m := \min \rset{k > n_{m-1} } {Y_k \ne Y_{k+1}}.
  \label{time-reparam}
\end{equation}
Then set $Y'_m := Y_{n_m}$. This reparametrization has the effect of
deleting those times when the $Y$-walker does not move. Manifestly,
$\{ Y'_m \}$ is the standard RW in the Manhattan lattice (at each
site, the two outgoing edges are equivalent, in terms of \pr).

By the result in the Appendix, $\{ Y'_m \}$ is recurrent. Therefore,
for infinitely many $n>0$, $Y_n = (0,0)$, that is, either $Z_n =
(0,0,0)$ or $Z_n = (1,1,1)$. If the second possibility happened all
the time, this would be equivalent to infinitely many failed attempts
at traveling from $(1,1,1)$ to $(0,0,0)$. Since, by the Markov
property and the hypothesis, these attempts are i.i.d.\ with positive
\pr, this situation can only occur with \pr\ zero.
\qed

It is clear from the previous proof that one can conjure up an
apparently similar \en\ with substantially different recurrence
properties. If $\{ \zeta_x \} \subset [0,1]$ as above, define, for $x
\in \Z^2_\mathrm{even}$,
\begin{equation}
  \begin{array}{rllll}
    (\mbox{if } d = \bfn, \bfs) && \omega_x(d, \bff) = \zeta_x, & 
    \omega_x(d, \bfl) = \omega_x(d, \bfr) = (1-\zeta_x)/2, 
    & \omega_x(d, \bfb) = 0; \\
    (\mbox{if } d = \bfw, \bfe) && \omega_x(d, \bff) = 0, & 
    \omega_x(d, \bfl) = \omega_x(d, \bfr) = 1/2, 
    & \omega_x(d, \bfb) = 0.
  \end{array}
  \label{hyp-4-1b-a}
\end{equation}
Also, for $x \in \Z^2_\mathrm{odd}$,
\begin{equation}
  \begin{array}{rllll}
    (\forall d \in \Delta) && \omega_x(d, \bff) = 0, & \omega_x(d,
    \bfl) = \omega_x(d, \bfr) = 1/2, & \omega_x(d, \bfb) = 0.
  \end{array}
  \label{hyp-4-1b-b}
\end{equation}
Notice that in this case the space average (\ref{avg-l-drift}) of the
local drift is zero. Yet the following is true:

\begin{proposition}
  Suppose that $\zeta_x > 0$ for at least one $x \in
  \Z^2_\mathrm{even}$. Then the recurrence of the above-defined \ra\
  walk is determined by its initial condition $\xi = (x_0,
  d_0)$. Namely, if
  \begin{displaymath}
    x_0 \in \Z^2_\mathrm{even},\ d_0 = \bfe, \bfw \qquad \mbox{or}
    \qquad x_0 \in \Z^2_\mathrm{odd},\ d_0 = \bfn, \bfs
  \end{displaymath}
  the walk is recurrent; otherwise it is non-recurrent. In either case
  the walk has zero asymptotic velocity.
  \label{prop-4-1b}
\end{proposition}

\proof Assumptions (\ref{hyp-4-1b-a})-(\ref{hyp-4-1b-b}) imply that
all the effective non-horizontal edges belong to $L_+$ (cf.\
Fig.~\ref{fig-rwlg-4}). Therefore, if the walk lands in $M_2$, it
cannot go back to $M_1$. Furthermore, by hypothesis, all walks
eventually land in $M_2$ (because---see the Appendix---the standard RW
in $M_1$ goes everywhere a.s., and thus will keep visiting sites $x$
with $\zeta_x > 0$ until it finally rises to $M_2$).

The proof is now immediate, as the assertion of Proposition
\ref{prop-4-1b} reads: If the initial condition is in $M_2$, the RW is
recurrent; if it is in $M_1$, the walk is not recurrent. The null
asymptotic velocity is a consequence of the properties of the stardard
RW in $M_2$.  
\qed

\subsection{Inhomogeneous Backward probability}

The examples presented above are not exceptional, and what can be done
with an inhomogeneous Forward probability can be done with an
inhomogeneous Backward probability as well.

For instance, take once again a collection $\{ \zeta_x \}$ of numbers
from $[0,1]$.  For $x \in \Z^2_\mathrm{even}$, set
\begin{equation}
  \begin{array}{rllll}
    (\mbox{if } d = \bfn) && \omega_x(d, \bff) = 0, & 
    \omega_x(d, \bfl) = \omega_x(d, \bfr) = (1-\zeta_x)/2, 
    & \omega_x(d, \bfb) = \zeta_x; \\
    (\mbox{if } d \ne \bfn) && \omega_x(d, \bff) = 0, & 
    \omega_x(d, \bfl) = \omega_x(d, \bfr) = 1/2, 
    & \omega_x(d, \bfb) = 0.
  \end{array}
  \label{hyp-4-2-a}
\end{equation}
For $x \in \Z^2_\mathrm{odd}$, set
\begin{equation}
  \begin{array}{rllll}
    (\mbox{if } d = \bfs) && \omega_x(d, \bff) = 0, & 
    \omega_x(d, \bfl) = \omega_x(d, \bfr) = (1-\zeta_x)/2, 
    & \omega_x(d, \bfb) = \zeta_x; \\
    (\mbox{if } d \ne \bfs) && \omega_x(d, \bff) = 0, & 
    \omega_x(d, \bfl) = \omega_x(d, \bfr) = 1/2, 
    & \omega_x(d, \bfb) = 0.
  \end{array}
  \label{hyp-4-2-b}
\end{equation}
As for the model of Proposition \ref{prop-4-1}, the local drift need
not be statistically balanced.

\begin{proposition}
  Suppose that, for at least one $x \in \Z^2_\mathrm{even}$, $\zeta_x
  > 0$ and, for at least one $y \in \Z^2_\mathrm{odd}$, $\zeta_y >
  0$. Suppose also that there is no $x \in \Z^2_\mathrm{even}$ such
  that $\zeta_x = \zeta_{x + \bfs} = 1$. Then the PRW in the \en\
  defined above is recurrent for all initial conditions.
  \label{prop-4-2}
\end{proposition}

\proof As anticipated, in this proof too we regard $\Gamma$ as
embedded in $\R^3$, with coordinates $(z^1, z^2, z^3)$.

By assumption, in this case, the effective edges of $L_+ \cup L_-$ are
those connecting $(z^1, z^2, 0)$ to $(z^1, z^2, 1)$ or viceversa, with
$(z^1, z^2) \in (2\Z)^2$. Define $\bar{\psi} : \Z^2 \times \{ 0,1 \}
\longrightarrow \Z^2$ as
\begin{equation}
  \bar{\psi}(z^1, z^2, z^3) := (z^1 + z^3, z^2 + z^3).
  \label{proj-psip}
\end{equation}
In analogy with (\ref{proj-psi}), $\bar{\psi}$ induces on $\Gamma$ a
translation of $M_2$ by the vector $(1,1,0)$ (again, this aligns the
orientations of $M_2$ and $M_1$), followed by a vertical projection
onto $M_1$. See Fig.~\ref{fig-rwlg-6}.

\fig{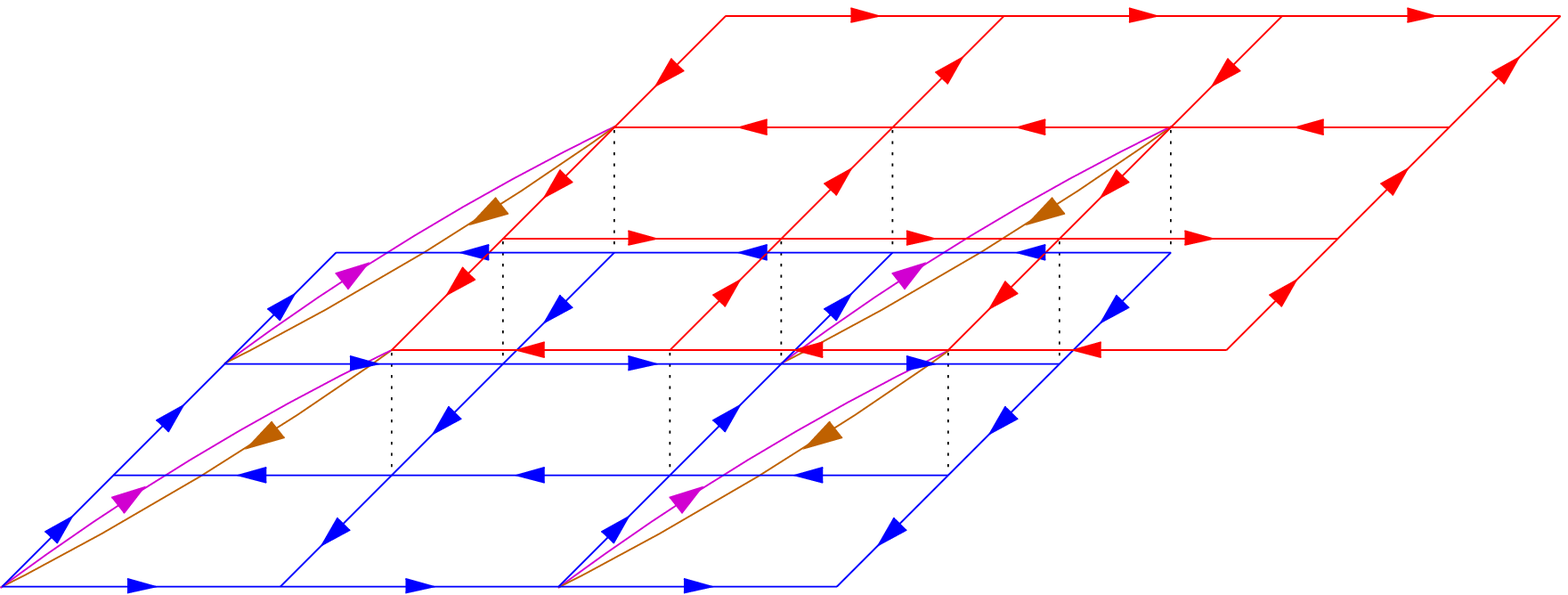}{4in}{Illustration of the proof of Proposition
\ref{prop-4-2}.}

The non-horizontal effective edges are projected into diagonals of
those unit squares of $\Z^2$ whose lower left corner belongs to
$(2\Z)^2$.  (These squares correspond, in the original problem, to the
sites of $\Z^2_\mathrm{even}$; compare Figs.~\ref{fig-rwlg-3} and
\ref{fig-rwlg-4}). These are the squares that, in the proof given
in the Appendix, are each collapsed to a single point (see
Fig.~\ref{fig-rwlg-7}).

This suggests the final step of our proof. For $\{ Z_n \}$ a
realization of our RW and $\ph$ as in (\ref{proj-phi}), set $Y_n :=
(\ph \circ \bar{\psi}) (Z_n)$. It is clear that $\{ Y_n \}$ has many
``dead times'', i.e., times $n$ at which $Y_n = Y_{n+1}$. The second
hypothesis of Proposition \ref{prop-4-2}, excluding the possibility
that the $Z$-walker remains trapped in a pair of twin edges of $B$,
implies that all sequences of consecutive dead times are finite with
\pr\ 1.  Hence we can eliminate them, as in the proof of Proposition
\ref{prop-4-1}, by introducing $Y'_m := Y_{n_m}$, where the
reparametrization is the same as (\ref{time-reparam}).

It is very easy to check that $\{ Y'_m \}$ is the standard RW in
2D. Its recurrence yields the recurrence of $\{ Z_n \}$ by the same
argument as in the the proof of Proposition \ref{prop-4-1}.  
\qed

It is trivial, at this point, to construct an \en\ with inhomogeneous,
statistically balanced, Backward probabilities, for which the
conclusions of Proposition \ref{prop-4-1b} hold. We will not bore the
reader with that.

\subsection{Random Left-Right environment}

In all the previous examples we have assumed that, at each site, the
\pr\ for the walker to go Left or Right was the same. We needed this
symmetry to somehow project the given RW into a standard RW in the
Mahnattan lattice, or similar.  With the next and last example, we go
beyond this assumption and show that the same techniques can be pushed
further, with the help of some results from the theory of ordinary
RWRE in $\Z^2$.

For $0 < \eps < 1/2$ and $x$ ranging in $\Z^2_\mathrm{even}$, let
$(\zeta_x, \zeta'_x)$ be i.i.d.\ \ra\ vectors from $[\eps, 1 -
\eps]^2$, w.r.t.\ to a \pr\ distribution whose details are of no
relevance. To each realization of $\{ \zeta_x, \zeta'_x \}$ associate
an \en\ $\omega \in \es$ as follows.  For $x \in \Z^2_\mathrm{even}$,
\begin{equation}
  (\forall d \in \Delta) \qquad \omega_x(d, \bfl) = \omega_x(d, \bfr) =
  1/2
  \label{hyp-4-3-a}
\end{equation}
and
\begin{equation}
  \begin{array}{ll}
    \omega_{x+\bfe} (\bfe, \bfn) = \omega_{x+\bfw} (\bfw, \bfs) =
    \zeta_x, & \omega_{x+\bfe} (\bfe, \bfs) = \omega_{x+\bfw} (\bfw,
    \bfn) = 1-\zeta_x; \\
    \omega_{x+\bfn} (\bfn, \bfe) = \omega_{x+\bfs} (\bfs, \bfw) =
    \zeta'_x, & \omega_{x+\bfn} (\bfn, \bfw) = \omega_{x+\bfs} (\bfs,
    \bfe) = 1-\zeta'_x.
  \end{array}
  \label{hyp-4-3-b}
\end{equation}
Finally, for all $x \in \Z^2$ and $d \in \Delta$,
\begin{equation}
  \omega_x(d, d) = \omega_x(d, -d) = 0.
  \label{hyp-4-3-c}
\end{equation}
The resulting RE is thus of the Left-Right type and is denoted, as
usual, by $(\es, \Sigma, \Pi)$.

\begin{proposition}
  The PRWRE defined above verifies the CLT and is recurrent.
  \label{prop-4-3}
\end{proposition}

\proof As discussed earlier, in a Left-Right \en\ the edges of $L_+
\cup L_-$ are not effective and the walk ranges entirely within $M_1$
or $M_2$, depending on the initial condition. Upon closer inspection
of (\ref{hyp-4-3-b}), we realize that the transition probabilities for
the edges of $M_1$ are determined solely by $\{ \zeta_x \}$, whereas
the probabilities for the edges of $M_2$ are governed by $\{ \zeta'_x
\}$.

Without loss of generality, assume that the initial condition is $z_0
= (0,0,0)$, which lies in $M_1$. This corresponds, in the original
picture, to $(x_0, d_0) = (0, \bfn) = \xi$.

We use the same technique as in the Appendix to map our PRWRE $\{ Z_n
\}$ into a known RWRE in $\Z^2$. For
\begin{equation}
  \bar{\ph} (z^1, z^2, z^3) := \left( \left[ \frac{z^1}2 \right] ,
  \left[ \frac{z^2}2 \right] \right),
  \label{proj-phip}
\end{equation}
set $Y_n := \bar{\ph} (Z_{2n})$. Let us derive the transition
probabilities that the laws of $\{ Z_n \}$ induce on $\{ Y_n \}$.
Call $y$ the element of $\Z^2$ associated to $(x,d) \in
\Z^2_\mathrm{even} \times \{ \bfn, \bfs \}$ via the correspondence
$(x,d) \mapsto z \mapsto y$.  With the visual aid of
Figs.~\ref{fig-rwlg-7} and \ref{fig-rwlg-3} we verify that
\begin{eqnarray}
  && \hspace*{-16pt} P_\xi^\omega \{ Y_{n+1} = y + \bfe \,|\, Y_n = y
  \} = \omega_x (d,\bfe) \: \omega_{x+\bfe} (\bfe,\bfn) \,=\,
  \zeta_x/2; \nonumber \\
  && \hspace*{-16pt} P_\xi^\omega \{ Y_{n+1} = y + \bfn \,|\, Y_n = y
  \} = \omega_x (d,\bfw) \: \omega_{x+\bfw} (\bfw,\bfn) \,=\,
  (1-\zeta_x)/2; \nonumber \\
  && \hspace*{-16pt} P_\xi^\omega \{ Y_{n+1} = y + \bfw \,|\, Y_n = y \} =
  \omega_x (d,\bfw) \: \omega_{x+\bfw} (\bfw,\bfs) \,=\, \zeta_x/2;
  \label{trans-pr-y} \\
  && \hspace*{-16pt} P_\xi^\omega \{ Y_{n+1} = y + \bfs \,|\, Y_n = y
  \} = \omega_x (d,\bfe) \: \omega_{x+\bfe} (\bfe,\bfs) \,=\,
  (1-\zeta_x)/2.  \nonumber
\end{eqnarray}
(We restrict our attention to $(x,d) \in \Z^2_\mathrm{even} \times \{
\bfn, \bfs \}$ because, with the chosen initial condition, $(X_{2n},
D_{2n}) \in \Z^2_\mathrm{even} \times \{ \bfn, \bfs \} \ \forall
n>0$.)

Equations (\ref{trans-pr-y}) show that the \en\ for the $Y$-walk is
balanced, in the sense of Lawler, and elliptic. It is known that,
$\Pi$-a.s., this ordinary RW in $\Z^2$ satisfies the CLT and is
recurrent \cite{la, z}. As for for the original walk, the CLT holds
because the time-reparametrization between $\{ Y_n \}$ and $\{ Z_n
\}$, or $\{ (X_n, D_n) \}$, is simply linear; and the recurrence is
given by the same argument used in the last few proofs.  \qed

\appendix

\section{Appendix: The Manhattan lattice}
\label{sec-manh}

The (two-dimensional) \emph{Manhattan lattice} is defined as the
oriented graph ${\mathbb M}$ whose vertex set is $\Z^2$ and whose edge
set is
\begin{eqnarray}
  E_{\mathbb M} &:=& \rset{ \, [ \, (x^1, x^2) \,,\, (x^1 +
  \epsilon(x^2), x^2 ) \, ] } { (x^1, x^2) \in \Z^2 } \, \cup 
  \nonumber \\
  && \rset{ \, [ \, (x^1, x^2) \,,\, (x^1, x^2 + \epsilon(x^1)) \, 
  ] } { (x^1, x^2) \in \Z^2 },
\end{eqnarray}
where $\epsilon(j) = 1$ or $-1$, depending on $j \in \Z$ being even or
odd, respectively. See Fig.~\ref{fig-rwlg-7}.

\fig{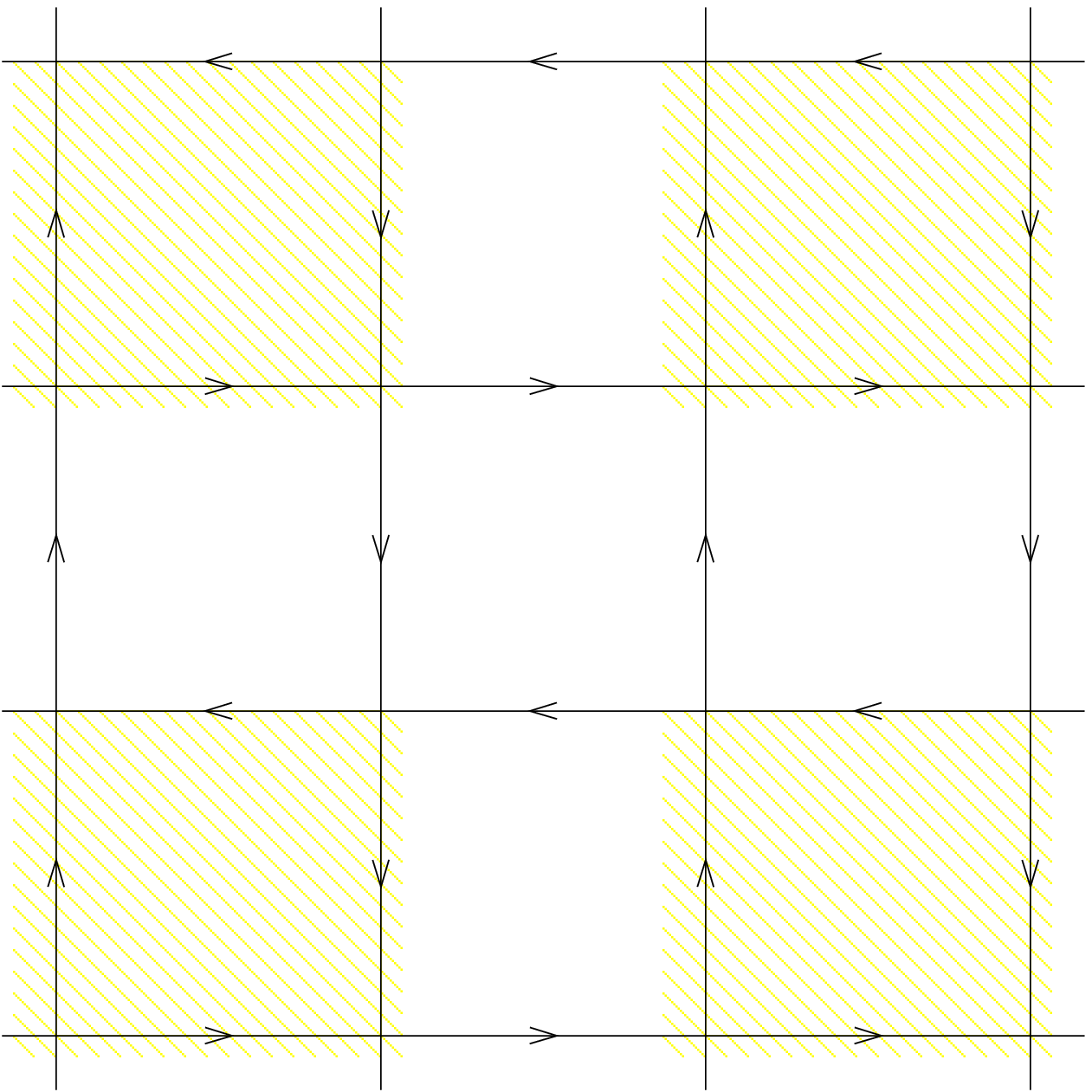}{3.2in}{The Manhattan lattice ${\mathbb M}$. If the
four vertices of each yellow square are identified and time is
accelerated by a factor 2, the standard \rw\ in ${\mathbb M}$ becomes
the standard \rw\ in $\Z^2$.}

The standard \rw\ in ${\mathbb M}$ is the one for which, at each node
$x \in \Z^2$, the \ra\ walker has \pr\ 1/2 of following either
outgoing edge based in $x$. We present a proof by Guillotin \cite{g}
that this RW has the same statistical features as the standard RW in
2D---in particular it satisfies the CLT and is recurrent. For $x =
(x^1, x^2) \in \Z^2$, define
\begin{equation}
  \ph(x) := \left( \left[ \frac{x^1}2 \right] , \left[ \frac{x^2}2
  \right] \right),
  \label{proj-phi}
\end{equation}
where $[j]$ is the largest integer less than or equal to $j \in
\Z$. If $\{ X_n \}_{n \in \N}$ is a realization of the RW in ${\mathbb
M}$, define $Y_n := \ph(X_{2n})$. Denote by $\Phi$ the transformation
$\{ X_n \} \mapsto \{ Y_n \}$, whose geometric meaning is explained in
the caption of Fig.~\ref{fig-rwlg-7}.

Fixing $X_0= (0,0)$ makes $\Phi$ invertible, and a little thought
convinces one that the endomorphism of \pr\ spaces that $\Phi$ induces
between the two RWs in $\Z^2$ is actually an automorphism. At any
rate, it is evident that $\{ Y_n \}$ is the stardard nearest-neighbor
RW in $\Z^2$, which immediately implies the CLT for $\{ X_n \}$.

As for the recurrence, one has that, a.s., $Y_n = (0,0)$ for
infinitely many $n$. For those $n$, either $X_{2n} = (0,0)$ or $X_{2n}
= (1,1)$ (\emph{tertium non datur}, as $2n$ is even).  By the Markov
property, the latter cannot happen all the time, because there is a
positive \pr\ of going from $(1,1)$ to $(0,0)$.

\end{document}